\def\draft{n}
\documentclass{amsart}
\usepackage[headings]{fullpage}
\usepackage{amssymb,epic,eepic,epsfig,amsmath,amsbsy,amsmath,bbm}
\usepackage{graphicx}
\usepackage{url}
\usepackage[bookmarks=true,%
    colorlinks=true,%
    linkcolor=blue,%
    citecolor=blue,%
    filecolor=blue,%
    menucolor=blue,%
    urlcolor=blue,%
    breaklinks=true]{hyperref}

\newtheorem{theorem}{Theorem}[section]
\newtheorem{proposition}{Proposition}[section]
\theoremstyle{definition}

\newtheorem{conjecture}[proposition]{Conjecture}

\def\printname#1{
        \if\draft y
                \smash{\makebox[0pt]{\hspace{-0.5in}
                        \raisebox{8pt}{\tt\tiny #1}}}
        \fi
}

\newcommand{\psdraw}[2]
         {\begin{array}{c} \hspace{-1.3mm}
        \raisebox{-4pt}{\epsfig{figure=draws/#1.eps,width=#2}}
       \hspace{-1.9mm}\end{array}}

\newlength{\standardunitlength}
\catcode`\@=11
\long\def\@makecaption#1#2{%
     \vskip 10pt

\setbox\@tempboxa\hbox{
       \small\sf{\bfcaptionfont #1. }\ignorespaces #2}%
     \ifdim \wd\@tempboxa >\captionwidth {%
         \rightskip=\@captionmargin\leftskip=\@captionmargin
         \unhbox\@tempboxa\par}%
       \else
         \hbox to\hsize{\hfil\box\@tempboxa\hfil}%
     \fi}
\font\bfcaptionfont=cmssbx10 scaled \magstephalf
\newdimen\@captionmargin\@captionmargin=2\parindent
\newdimen\captionwidth\captionwidth=\hsize
\catcode`\@=12

\def\lbl#1{\label{#1}\printname{#1}}

\def\BN{\mathbbm N}
\def\BZ{\mathbbm Z}
\def\BQ{\mathbbm Q}

\def\BC{\mathbbm C}

\def\D{\Delta}

\def\a{\alpha}

\def\la{\langle}
\def\ra{\rangle}

\def\e{\epsilon}

\def\d{\delta}
\def\b{\beta}

\def\fsl{\mathfrak{sl}}

\newcommand{\qbinom}[2]{\genfrac{[}{]}{0pt}{}{#1}{#2}}
\def\Tet{\mathrm{Tet}}
\def\U{\mathrm{U}}

\def\vol{\mathrm{vol}}

\begin{document}


\title[The non-commutative $A$-polynomial of $(-2,3,n)$ pretzel knots]{
The non-commutative $A$-polynomial of $(-2,3,n)$ pretzel knots}
\author{Stavros Garoufalidis}
\address{School of Mathematics \\
         Georgia Institute of Technology \\
         Atlanta, GA 30332-0160, USA \newline 
         {\tt \url{http://www.math.gatech.edu/~stavros}}}
\email{stavros@math.gatech.edu}
\author{Christoph Koutschan}
\address{Research Institute for Symbolic Computation \\
         Johannes Kepler University \\
         Altenbergerstrasse 69 \\
         A-4040 Linz, Austria 
         \tt{\url{http://www.risc.jku.at/home/ckoutsch}}}
\email{Koutschan@risc.uni-linz.ac.at}
\thanks{C.K. was supported by grants DMS-0070567 of US National Science 
Foundation and FWF P20162-N18 of the Austrian Science Fund.
S.G. was supported in part by grant DMS-0805078 of the US National Science 
Foundation.
\\
\newline
1991 {\em Mathematics Classification.} Primary 57N10. Secondary 57M25.
\newline
{\em Key words and phrases: colored Jones polynomial, knots, pretzel
knots, non-commutative $A$-polynomial, $q$-holonomic sequences,
recursion ideal, Quantum Topology, Volume Conjecture, Kashaev invariant.
}
}

\date{December 7, 2011}


\begin{abstract}
We study $q$-holonomic sequences that arise as the colored Jones
polynomial of knots in 3-space.  The minimal-order recurrence for such
a sequence is called the (non-commutative) $A$-polynomial of a knot.
Using the \emph{method of guessing},
we obtain this polynomial explicitly for the $K_p=(-2,3,3+2p)$ pretzel knots 
for $p=-5,\dots,5$. This is a particularly interesting family since the
pairs $(K_p,-K_{-p})$ are geometrically similar (in particular, scissors 
congruent) with similar character varieties. Our computation of the
non-commutative $A$-polynomial 
(a) complements the computation of the $A$-polynomial of the pretzel knots
done by the first author and Mattman, (b) supports the AJ Conjecture for
knots with reducible $A$-polynomial and (c) numerically computes 
the Kashaev invariant of pretzel knots in linear time. In a later
publication, we will use the numerical computation of the Kashaev invariant
to numerically verify the Volume Conjecture for the above mentioned 
pretzel knots.
\end{abstract}

\maketitle


\section{The colored Jones polynomial: a $q$-holonomic 
sequence of natural origin}
\lbl{sec.intro}

\subsection{Introduction}
\lbl{sub.introduction}

The colored Jones polynomial of a knot $K$ in 3-space is a $q$-holonomic
sequence of Laurent polynomials of natural origin in Quantum Topology
\cite{GL1}. As a canonical recursion relation for this sequence we choose the
one with minimal order; this is the so-called non-commutative $A$-polynomial 
of a knot \cite{Ga1}. Using the computational \emph{method of guessing} with
undetermined coefficients
\cite{Ka1,Ka2} combined with a carefully chosen exponent set of monomials
(given by a translate of the Newton polygon of the $A$-polynomial)
we compute very plausible candidates for the non-commutative $A$-polynomial 
of the $(-2,3,3+2p)$ pretzel knot family for $p=-5,\dots,5$. Our computation
of the non-commutative $A$-polynomial 
\begin{itemize}
\item[(a)]
complement the computation of the $A$-polynomial of the pretzel knots
\cite{GM}, 
\item[(b)]
support the AJ Conjecture of \cite{Ga1} (see also \cite{Ge}) for knots with 
reducible $A$-polynomial, and
\item[(c)] 
give an efficient linear time algorithm for computing numerically the Kashaev 
invariant of the pretzel knots (with a fixed accuracy).
\end{itemize}
In \cite{GZ}, we use the latter algorithm to numerically verify the 
volume conjecture of Kashaev \cite{Ks,MM} for the above 
mentioned pretzel knots. 

For an introduction to the polynomial invariants of knots that
originate in Quantum Topology \cite{Jo,Tu1,Tu2} and the book \cite{Ja}
where all the details of the quantum group theory can be found. For
up-to-date computer calculations of several polynomial invariants of knots, 
see \cite{B-N}. For an introduction to $q$-holonomic sequences see
\cite{Z,WZ,PWZ}. For the appearance of $q$-holonomic sequences in Quantum
Topology, see \cite{GL1,GS1,GS2} and also \cite{GK}.

\subsection{Fusion and the colored Jones polynomial of pretzel knots}
\lbl{sub.pretzel}

Consider the 1-parameter family of {\em pretzel knots} $K_p=(-2,3,3+2p)$
for an integer $p$
$$
\psdraw{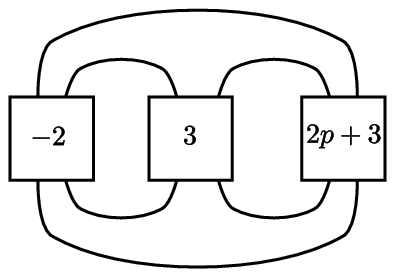}{2.0in}
$$
where an integer $m$ inside a box indicates the number $|m|$ of
half-twists, right-handed (if $m>0$) or left-handed (if $m<0$), according
to the following figure:
$$
\psdraw{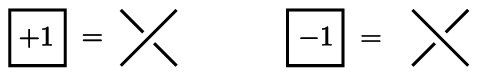}{2.5in}
$$
The pretzel knots $K_p$ are interesting from many points of view, discussed
in detail in \cite{CGLS,GM,Ga5,GZ}:
\begin{itemize}
\item
In {\em hyperbolic geometry}, $K_p$ is the torus knot $5_1, 8_{19}$ and
$10_{124}$ when $p=-1,0,1$, and $K_p$ is a hyperbolic knot when $p \neq -1,0,1$.
\item
The pairs $(K_p,-K_{-p})$ (where $-K$ denotes the mirror of $K$) are
{\em geometrically similar} for $p \geq 2$. In particular, their complements 
are scissors congruent, with equal volume, and with Chern-Simons 
invariants differing by torsion \cite{Ga5}.
\item
The knots $K_p$ appear in the study of {\em exceptional Dehn surgery}
\cite{CGLS}.
\item
In {\em Quantum Topology}, the knots $K_p$ have different Jones polynomial
and different Kashaev invariants, which numerically verify the Volume
Conjecture \cite{GZ}.
\end{itemize}
Let $J_{K,n}(q)$ denote the colored {\em colored Jones polynomial} 
of a knot $K$ colored by the $n$-dimensional irreducible representation of
$\fsl_2$, framed by zero and normalized to be $1$ at the unknot
\cite{Tu1,Tu2}. So, $J_{K,1}(q)=1$ for all knots and $J_{K,2}(q)$
is the Jones polynomial of $K$ \cite{Jo}. 
Our starting point is an explicit formula for the colored Jones polynomial 
$J_{p,n}(q)$ of $K_p$. This comes from a theorem of \cite{Ga2} which has two 
parts. The first part 
identifies the pretzel knots $K_p$ with members of a 2-parameter family of 
2-{\em fusion knots} $K(m_1,m_2)$ for integers $m_1$ and $m_2$, drawn here
$$
\psdraw{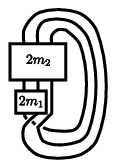}{1in}
$$
and discussed in detail in \cite{Ga2}. 
The second part gives an explicit formula for the colored
Jones polynomial of $K(m_1,m_2)$. To state it, we need to recall some notation.
The {\em quantum integer} $[n]$ and the {\em quantum factorial} 
$[n]!$ of a natural number $n$ are defined by
$$
[n]=\frac{q^{n/2}-q^{-n/2}}{q^{1/2}-q^{-1/2}},
\qquad
[n]!=\prod_{k=1}^n [k]!
$$
with the convention that $[0]!=1$. Let 
$$
\qbinom{a}{a_1, a_2, \dots, a_r}=\frac{[a]!}{[a_1]! \dots [a_r]!}
$$
denote the $q$-multinomial coefficient of natural numbers $a_i$ such that
$a_1+ \dots +a_r=a$. We say that a triple $(a,b,c)$ of natural
numbers is {\em admissible} if $a+b+c$ is even and the triangle inequalities
hold. In the formulas below, we use the following basic trivalent graphs
$\U,\Theta,\Tet$ colored by one, three and six natural numbers (one
in each edge of the corresponding graph) such that
the colors at every vertex form an admissible triple.

\begin{figure}[htpb]
$$
\psdraw{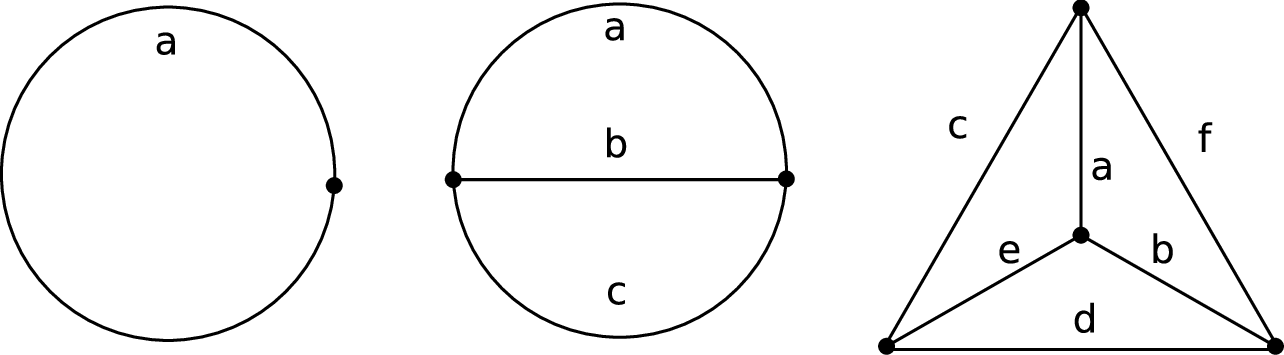}{3in}
$$
\caption{The $\U$, $\Theta$ and $\Tet$ graphs colored by an admissible 
coloring.}
\lbl{eq.3graphs}
\end{figure}
If a coloring of a graph is not admissible, its evaluation vanishes.
When the colorings of the graphs in Figure \eqref{eq.3graphs} are admissible,
their evaluations can be computed by the following functions.

\begin{eqnarray*}
\mu(a)&=& (-1)^a q^{\frac{a(a+2)}{4}} \\
\nu(c,a,b)&=& (-1)^{\frac{a+b-c}{2}} q^{\frac{-a(a+2)-b(b+2)+c(c+2)}{8}} \\
\U(a)&=&(-1)^a[a+1] \\
\Theta(a,b,c)&=& (-1)^{\frac{a+b+c}{2}}[\frac{a+b+c}{2}+1]
\qbinom{\frac{a+b+c}{2}}{\frac{-a+b+c}{2}, \frac{a-b+c}{2}, \frac{a+b-c}{2}}
\\
\Tet(a,b,c,d,e,f)&=& 
\sum_{k = \max T_i}^{\min S_j} (-1)^k [k+1]
\qbinom{k}{S_1-k , S_2-k , S_3-k , k- T_1 , k- T_2 , k- T_3 , k- T_4}
\end{eqnarray*}
where
\begin{equation}
\lbl{eq.Sj}
S_1 = \frac{1}{2}(a+d+b+c)\qquad S_2 = \frac{1}{2}(a+d+e+f) 
\qquad S_3 = \frac{1}{2}(b+c+e+f)
\end{equation}
\begin{equation}
\lbl{eq.Ti}
T_1 = \frac{1}{2}(a+b+e) \qquad T_2 = \frac{1}{2}(a+c+f)
\qquad T_3 = \frac{1}{2}(c+d+e) \qquad T_4 = \frac{1}{2}(b+d+f).
\end{equation}
An assembly of the five building blocks can compute the colored Jones
function of any knot. 
Consider the rational convex plane polygon $P$ 
with vertices $\{(0,0),(1/2,-1/2),(1,0),(1,1)\}$ in $\BQ^2$: 
$$
\psdraw{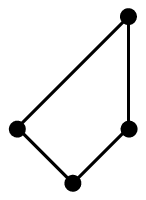}{0.5in}
$$

\begin{theorem}
\lbl{thm.fusion}\cite{Ga2}
\rm{(a)}
For every integer $p$, we have $K_p=K(p,1)$.
\newline
\rm{(b)} For every $m_1,m_2 \in \BZ$ and $n \in \BN$, we have:
\begin{eqnarray}
\lbl{eq.cjk}
J_{K(m_1,m_2),n+1}(1/q) &=& \frac{\mu(n)^{-w(m_1,m_2)}}{\U(n)} \sum_{(k_1,k_2) \in nP \cap \BZ^2}
\nu(2k_1,n,n)^{2 m_1+2 m_2} \nu(n+2 k_2,2k_1,n)^{2m_2+1} \\
& & \notag
\cdot \frac{\U(2k_1)\U(n+2k_2)}{\Theta(n,n,2k_1) \Theta(n,2k_1,n+2k_2)}
\Tet(n,2k_1,2k_1,n,n,n+2k_2) 
\end{eqnarray}
where $P$ is as above and 
the {\em writhe} of $K(m_1,m_2)$ is given by $w(m_1,m_2)=2m_1+6m_2+2$.
\end{theorem}

\subsection{Our results}
\lbl{sub.results}

Recall that a $q$-{\em holonomic sequence} $(f_n(q))$ for $ n \in \BN$ 
is a sequence (typically of rational functions $f_n(q) \in \BQ(q)$ in 
one variable $q$)
which satisfies a {\em linear recursion relation}:
\begin{equation}
\lbl{eq.recf}
a_d(q^n,q) f_{n+d}(q) + \dots + a_0(q^n,q) f_n(q) = b(q^n,q)
\end{equation}
for all $n \in \BN$, where $a_j(u,v) \in \BQ[u,v]$ for all $j=0,\dots,d$
and $b(u,v) \in \BQ[u,v]$ \cite{Z}. As is custom, one can phrase 
Equation \eqref{eq.recf} in operator form, by considering the operators
$M$ and $L$ that act on a sequence $(f_n(q))$ by:
\begin{align*}
(L f)_{n}(q)& = f_{n+1}(q) & (M f)_{n}(q)& = q^{n} f_{n}(q)
\end{align*}
It is easy to see that the operators $M$ and $L$ satisfy the $q$-commutation
relation:
$$
LM=qML
$$
Thus, we can write Equation \eqref{eq.recf} in the form:
\begin{equation}
\lbl{eq.recfP}
P f=b, \qquad P=\sum_{j=0}^d a_j(M,q)L^j, \qquad b=b(q^n,q)
\end{equation}
We will call a $q$-holonomic bi-infinite sequence $f_n(q)$ {\em palindromic}
if either $f_n(q)=f_{-n}(q)$ for all integers $n$, or 
$f_n(q)=-f_{-n}(q)$ for all integers $n$.
Given a palindromic sequence $f_n(q)$, 
we will call a recursion relation \eqref{eq.recf}
(and the corresponding operator $(P,b)$) {\em palindromic} if Equation
\eqref{eq.recf} holds for all integers $n$. 

With our normalizations, the colored Jones polynomial $J_{K,n}(q)$ of a knot, 
defined for $n \geq 1$, extends to a palindromic sequence defined by
$J_{K,n}(q)=J_{K,-n}(q)$ for $n <0$ and $J_{K,0}(q)=1$.

Let $A_p(M,L) \in \BQ[M^2,L]$ denote the {\em $A$-polynomial} of the pretzel 
knot $K_p$, given in \cite{GM}. Let $\e_p(M) \in \BQ[M]$ denote the 
$M$-factors given in Appendix \ref{sec.Mfactors}.
Let $\D_p(t) \in \BZ[t^{\pm 1}]$ denote the {\em Alexander polynomial} of 
$K_p$; \cite{Kf}. $\D_p$ satisfies a linear recursion 
relation:
\begin{equation}
\lbl{eq.recalex}
\D_{p+2}-(t+t^{-1}) \D_{p+1}+\D_p=0
\end{equation}
for all $p \in \BZ$, with initial conditions
$$
\D_0=\frac{1}{t^3}-\frac{1}{t^2}+1-t^2+t^3,
\qquad
\D_1=\frac{1}{t^4}-\frac{1}{t^3}+\frac{1}{t}-1+t-t^3+t^4.
$$

\begin{theorem}
\lbl{thm.1}
\rm{(a)}
Consider the operators $(A_p(M,L,q),b_p(M,q))$ of the appendix 
for $p=-5,\dots,5$. Then, we have:
$$
A_p(M,L,1)=A_p(M^{1/2},L) \e_p(M)
$$
in accordance with the AJ Conjecture \cite{Ga1}.
\newline
\rm{(b)}
$(A_p(M,L,q),b_p(M,q))$ is palindromic.
\newline
\rm{(c)}
We also have:
\begin{equation}
\lbl{eq.loop}
\frac{A_p(M,1,1)}{b_p(M,1)}=\D_p(M)
\end{equation}
for $p \neq -3$. When $p=-3$ we have $A_{-3}(M,1,1)=b_{-3}(M,1)=0$ and
\begin{equation}
\lbl{eq.loop-3}
\frac{A_{-3,q}(M,1,1)}{\D_{-3}(M)}-A_{-3,L}(M,1,1)M \frac{\D_{-3}'(M)}{\D_{-3}(M)^2}=
b_{-3,q}(M,1)
\end{equation}
where primes indicate partial derivatives. This is
in accordance with the loop expansion of the colored Jones 
polynomial \cite{Ga4}. For a definition of the loop expansion \cite{Ro}.
\end{theorem}

\begin{conjecture}
\lbl{conj.1}
We conjecture that: 
\begin{equation}
\lbl{eq.conjrec}
A_p(M,L,q) J_{p,n}(q) = b_p(q^n,q)
\end{equation}
for $p=-5,\dots, 5$ and all $n \in \BZ$. 
\end{conjecture}

\section{Consistency checks}
\lbl{sec.remarks}

Three consistency checks of Conjecture \ref{conj.1} were already mentioned
in Theorem \ref{thm.1}. In this section we discuss four independent
consistency checks regarding Conjecture \ref{conj.1}.

\subsection{Consistency with the height}
\lbl{sub.height}

In this section we discuss the height of Equation \eqref{eq.conjrec}.
For fixed integers $p$ and natural numbers $n$, both sides of Equation
\eqref{eq.conjrec} are Laurent polynomials in $q$ with integer
coefficients.

In general, the minimum and maximum degree of a $q$-holonomic sequence 
of Laurent
polynomials is a quadratic quasi-polynomial \cite{Ga7}.
In \cite{Ga2} the first author studied the minimum and maximum
degree of the colored Jones polynomial of the 2-fusion knots $K(m_1,m_2)$.
For the case of pretzel knots $K_p$, the maximum degree of $J_{p,n}(q)$
is a quadratic quasi-polynomial of $n$, and the minimum degree is a 
linear function of $n$. It follows that the terms in the left hand side of 
Equation \eqref{eq.conjrec} are polynomials of $q$ of minimum degree a linear
function of $n$ and maximum degree a quasi-polynomial quadratic function of 
$n$. 
Explicitly, for the case of $K_2=K(2,1)$, it was shown in \cite{Ga2,Ga6}
that $J_{2,n}(q)$ is a polynomial of $q$ of minimum degree $\d^*(n)$ and 
maximum degree $\d(n)$ given by:
\begin{eqnarray*}
\d(n) &=& \left[\frac{37}{8} n^2 + \frac{3}{4} n -\frac{31}{8}\right]=
\frac{37}{8} n^2 + \frac{3}{4} n -\frac{31}{8} + \e(n),
\\
\d^*(n) &=& 5(n-1)
\end{eqnarray*}
where $\e(n)$ is a periodic sequence of period $4$ given by $1/8,0,1/8,1/2$
if $n\equiv 0,1,2,3 \bmod 4$ respectively. 
Keep in mind that $J_{K,n}(q)$ denotes the colored Jones polynomial
of $K$ colored by the $n$-dimensional irreducible representation of
$\mathfrak{sl}_2$. It follows that for $p=2$,
the left hand side of Equation \eqref{eq.conjrec} is a polynomial in $q$ 
of minimum degree $60n+O(1)$ and maximum degree $37/8 n^2+9/4n+O(1)$.
We computed $J_{2,n}(q)$ explicitly for $1\leq n\leq 70$ using
Theorem \ref{thm.fusion}. For instance, $J_{2,70}(q)$ is a polynomial
of 
\begin{itemize}
\item
minimum (resp. maximum) exponent $345$ (resp. $22606$), 
\item
maximum (in absolute value) coefficient $14287764770955$ and 
\item
sum of absolute values of its coefficients $28587411833908277$.
\end{itemize} 
It follows that that Equation \eqref{eq.conjrec} for $n=70$ (which
actually holds, by an explicit computation) involves the matching of about 
$22500$ many powers of $q$ with coefficients $14$ digit integers.
One can compare this to the modest size of $(A_2(M,L,q),b_2(M,q))$ given in 
Appendix \ref{sec.K2}.

Of course, 
one can come up with operators that satisfy parts (a), (b) and (c)
of Theorem \ref{thm.1} and Equation \eqref{eq.conjrec} for
a fixed natural number $n$ (such as $n=70$). This is simply a problem of linear
algebra with more unknowns than coefficients which 
in fact has infinitely many solutions.
On the other hand, the operators given in the appendix are of small height, 
given the height of the input, as was illustrated above. 
Note also that Theorem \ref{thm.1} can be proven rigorously if we knew
a priori bounds for the $M,L$ degrees of the operators involved, and
if we were able to compute enough values of the colored Jones polynomials,
using the formula of Theorem \ref{thm.fusion}.

\subsection{Consistency with the loop expansion of the colored Jones 
polynomial}
\lbl{sub.loop}

This section concerns the consistency of Conjecture \ref{conj.1}
with the loop expansion of the colored Jones polynomial of $K_p$.
The latter was introduced by Rozansky in \cite{Ro}, and has the form:
\begin{equation}
\lbl{eq.loopK}
J_{K,n}(q)=\sum_{k=0}^\infty \frac{P_{K,k}(q^n)}{\D_K(q^n)^{2k+1}}(q-1)^k \in \BQ[[q-1]]
\end{equation}
where $P_k(t) \in \BZ[t^{\pm 1}]$ are Laurent polynomials with
$P_{K,0}=1$ and $\D_K(t) \in \BZ[t^{\pm 1}]$ is the Alexander polynomial of $K$.
With some effort, one can compute the $k$-loop polynomials $P_{K,k}(M)$ 
of a knot for various small values of $k$. The loop expansion given by
Equation 
\eqref{eq.loopK} {\em highly constrains} the coefficients of 
$(A_K(M,L,q),b_K(M,q))$, as was discussed in detail in \cite{Ga4}. The simplest
constraint for the knots $K_p$ is given in Equation \eqref{eq.loop}, which 
in fact is the first of a hierarchy of constraints. Each such constraint
gives a consistency check for Conjecture \ref{conj.1}, and offers a practical
way to compute the loop expansion of the knots $K_p$. Consistency with the
higher loop constraints have also been checked. We plan to
discuss their details in a forthcoming publication.

\subsection{Consistency with the AJ Conjecture}
\lbl{sub.AJ}

This section concerns the AJ Conjecture of \cite{Ga1} for the knots 
$K_{\pm 3}$ with reducible $A$-polynomial.
The $A$-polynomial $A_p(M,L)$ of the pretzel knots $K_p$ was computed
in \cite{GM}. In \cite{Ma} it was shown that $A_p(M,L)$ is irreducible 
if $3$ does not divide $p$ and otherwise it is the product of two
irreducible factors when $3$ divides $p$. Explicitly, we have:

\vspace{0.1in}
\begin{math}
A_{-3}(M,L)=
-(-1+L) (L^3-L^4-5 L^3 M^2+L^4 M^2-2 L^2 M^4-2 L^4 M^4-L M^6-4 L^2 M^6+3 L^3 M^6+2 L^4 M^6-M^8-5 L M^8-3 L^3 M^8+L^4 M^8+L M^{10}-3 L^2 M^{10}-5 L^4 M^{10}-L^5 M^{10}+2 L M^{12}+3 L^2 M^{12}-4 L^3 M^{12}-L^4 M^{12}-2 L M^{14}-2 L^3 M^{14}+L M^{16}-5 L^2 M^{16}-L M^{18}+L^2 M^{18})
\end{math}

\vspace{0.1in}
\begin{math}
A_3(M,L)=
(-1+L M^{24}) (-1+L M^{16}-L M^{18}+2 L M^{20}-5 L M^{22}+L M^{24}+5 L^2 M^{40}-4 L^2 M^{42}+L^2 M^{46}+L^3 M^{62}+3 L^3 M^{66}+2 L^3 M^{68}-2 L^4 M^{84}-3 L^4 M^{86}+3 L^4 M^{88}+2 L^4 M^{90}-2 L^5 M^{106}-3 L^5 M^{108}-L^5 M^{112}-L^6 M^{128}+4 L^6 M^{132}-5 L^6 M^{134}-L^7 M^{150}+5 L^7 M^{152}-2 L^7 M^{154}+L^7 M^{156}-L^7 M^{158}+L^8 M^{174})
\end{math}

\vspace{0.1in}
Theorem \ref{thm.1} matches exactly with the above values of the 
$A$-polynomials.
The case of the knot $K_{-3}$ is particularly interesting, since its
$\mathrm{SL}_2(\BC)$ character variety has {\em three} components: 
the geometric one, the abelian $L-1$ component, and an additional 
$L-1$ component of 
nonabelian representations. Theorem \ref{thm.1} and Conjecture \ref{conj.1}
support the idea that the AJ Conjecture captures the multiplicity of the
various components of the character variety.

\subsection{Consistency with the Volume Conjecture}
\lbl{sub.vc}

This section concerns the consistency of Conjecture \ref{conj.1} with
the computation of the Kashaev invariant of the $K_p$ knots.
The $N$-th {\em Kashaev invariant} $\la K \ra_N$ of a knot $K$ is defined
by \cite{Ks,MM}:
\begin{equation}
\lbl{eq.kashaev}
\la K \ra_N= J_{K,N}(e^{2 \pi i/N})
\end{equation}
The Volume Conjecture of Kashaev states that if $K$ is a hyperbolic knot,
then
\begin{equation}
\lbl{eq.vc}
\lim_{N \to \infty} \frac{|\la K \ra_N|}{N}=\frac{\text{vol}(K)}{2\pi}
\end{equation}
where $\text{vol}(K)$ is the volume of the hyperbolic knot $K$.
Since we are specializing to a root of unity, we might as well consider
the remainder $\tau_{K,N}(q)$ of $J_{K,N-1}(q)$ by 
the $N$-th cyclotomic polynomial $\Phi_N(q)$. In \cite{GZ}, it was shown that
given a recursion relation for $J_{K,N}(q)$, there is a
linear time algorithm to numerically compute $\la K \ra_N$.  
Using the guessed recursion relation for $K_2$, we compute $\tau_{K_2,N}(q)$
for $N=1,\dots,1000$. Here is a sample computation.

{\small

\vspace{0.1in}
\begin{math}
\tau_{K_2,100}(q)=\\
-1420771679897311607360-1402034476570732425908 q
-1377764083694494707679 q^2\\
-1348056285420017550322 q^3-1313028324854995190830 q^4
-1272818441358081463973 q^5\\-1227585324968178744317 q^6
-1177507490130630983388 q^7-1122782571182284245313 q^8\\
-1063626542375688303231 q^9+420498814366636734411 q^{10}
+469062907903390306537 q^{11}\\+515775824438145014436 q^{12}
+560453209429428890901 q^{13}+602918741648741441924 q^{14}\\
+643004829043136905736 q^{15}+680553270138355921566 q^{16}
+715415878390451489264 q^{17}\\+747455067013913965248 q^{18}
+776544391967778302155 q^{19}-618202628922511743188 q^{20}\\
-576608139973286430388 q^{21}-532738042123286363977 q^{22}
-486765470606610517117 q^{23}\\-438871858158259827294 q^{24}
-389246218987652812332 q^{25}-338084402821172432280 q^{26}\\
-285588321971646221647 q^{27}-231965154488540570326 q^{28}
-177426526516296620808 q^{29}\\+1298584002796105745794 q^{30}
+1335567867823634101034 q^{31}+1367280856639633305993 q^{32}\\
+1393597812566394292363 q^{33}+1414414874600710903331 q^{34}
+1429649887309469255114 q^{35}\\+1439242725058651352936 q^{36}
+1443155529298983637839 q^{37}+1441372857979981026638 q^{38}\\
+1433901746491878528487 q^{39}
\end{math}
}

Let 
$$
a_N=2 \pi \frac{\log|\la K_2 \ra_{N}|}{N}
$$
Since $\la K_2 \ra_N=\tau_{K_2,N}(e^{2 \pi i/N})$, the above expression 
gives the numerical value:
$$
a_{100}= 3.22309\dots
$$
which is a rather poor approximation of the volume 
$\vol(K_2)=2.8281220883307827\dots$ of $K_2$. On the other hand,
for $N=990,\dots,1000$ we have:

\begin{center}
{\small
\begin{tabular}{|c|c|c|c|c|c|c|c|c|c|c|c|}\hline
$N$ & $990$ & $991$ & $992$ & $993$ & $994$ & $995$ & $996$ & $997$ & $998$ & $999$ & $1000$ \\\hline
$a_N$ & $2.88981$ & $2.88976$ & $2.88971$ & $2.88965$ & $2.8896$ & $2.88955$ & $2.8895$ & $2.88944$ & $2.88939$ & $2.88934$ & $2.88929$ \\ \hline
\end{tabular}
}
\end{center}
The above data plots as follows:
$$
\psdraw{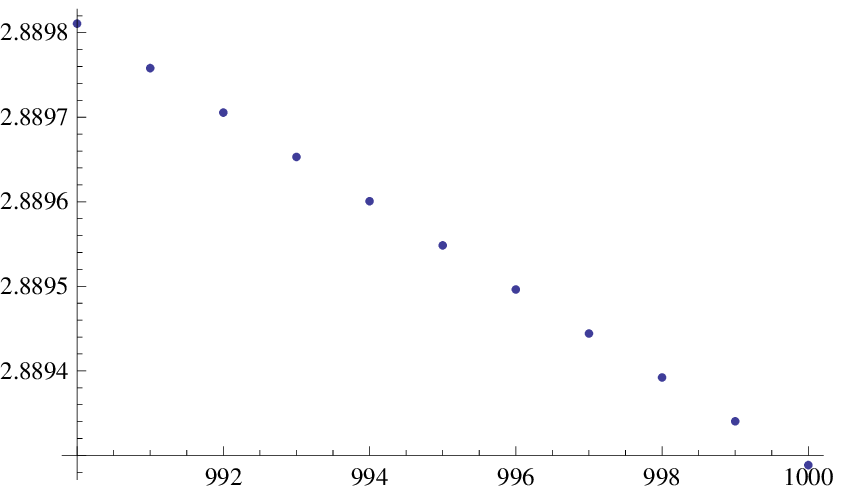}{3in}
$$
and numerically fits the following curve:
$$
2.82813 + 9.41764 \frac{\log(n)}{n} - 3.89193 \frac{1}{n} 
$$
which is a 4-digit approximation to the volume. In \cite{GZ} a more
precise approximation to the volume and its correction is given.

\section{The computation}
\lbl{sec.thm1}

We use the computer to guess the recurrences
for the colored Jones polynomials~$J_{p,n}(q)$ \cite{Ka1,Ko1,Ko2,GK}.
The term \emph{guessing} here refers to the method of making an ansatz 
with undetermined coefficients. The first values of
the sequence~$J_{p,n}(q)$ can be computed explicitly using Equation~\eqref{eq.cjk}
(see Table~\ref{tab.cjk0} for an example).  These values are then
plugged into an ansatz with undetermined coefficients in order to
produce an overdetermined linear system of equations. The more
equations are used (note that at least as many equations as unknowns
are needed to obtain a reliable result), the higher is the certainty
that the result is the recurrence for the sequence.
Alternatively, we can verify that the recurrence is 
satisfied for some values of the sequence that were not used for 
guessing, gaining further confidence into the result.

\begin{table}[htpb]
\begin{tabular}{|l|p{0.9\textwidth}|}\hline
$n$ & $J_{0,n}(q)$ \\ \hline\hline
$1$\rule{0em}{1em} & $1$\\ \hline
$2$\rule{0em}{1em}  & $-q^8+q^5+q^3$\\ \hline
$3$\rule{0em}{1em}  & $q^{23}-q^{22}+q^{20}-q^{19}-q^{16}-q^{13}+q^{12}+q^9+q^6$\\ \hline
$4$\rule{0em}{1em}  & $-q^{43}+q^{41}+q^{40}-q^{39}+q^{37}-q^{35}+q^{33}-q^{31}+q^{29}-q^{27}-q^{26}+
      q^{25}-q^{23}-q^{22}+q^{21}-q^{19}+q^{17}+q^{13}+q^9$\\ \hline
$5$\rule{0em}{1em}  & $q^{70}-q^{69}+q^{65}-2q^{64}+q^{60}-q^{59}+q^{57}+q^{55}-q^{54}+q^{52}-q^{49}+
      q^{47}-q^{44}+q^{42}-q^{39}+q^{37}-q^{35}-q^{34}+q^{32}-q^{30}-q^{29}+q^{27}-
      q^{25}+q^{22}+q^{17}+q^{12}$ \\ \hline
\end{tabular}
\caption{The first elements of the colored Jones polynomial of the (-2,3,3) pretzel knot~$K_0$}
\label{tab.cjk0}
\end{table}

\begin{table}[htpb]
\begin{tabular}{|r|c|c|c|}\hline
$p$\rule{0em}{1em}  & $d(J_{p,10}(q))$ & $d(J_{p,20}(q))$ & $d(J_{p,30}(q))$\\ \hline
-5\rule{0em}{1em} & 453 & 1919 & 4400 \\
-4 & 363 & 1546 & 3549 \\
-3 & 282 & 1197 & 2735 \\
-2 & 225 & 950 & 2175 \\
-1 & 225 & 950 & 2175 \\
0 & 265 & 1130 & 2595 \\
1 & 330 & 1410 & 3240 \\
2 & 406 & 1736 & 3991 \\
3 & 491 & 2098 & 4821 \\
4 & 579 & 2469 & 5671 \\
5 & 667 & 2843 & 6529 \\ \hline

\end{tabular}
\caption{Size of the colored Jones polynomial at $n=10,20,30$ for the pretzel knot family,
where $d(p)=d_1+d_2$ for a Laurent polynomial $p=\sum_{i=-d_1}^{d_2}c_iq^i$
with $c_{-d_1}\neq0$ and $c_{d_2}\neq0$}
\label{tab.Jdeg}
\end{table}

For guessing $q$-difference equations, there are two different choices 
for the ansatz and the nature of its coefficients. 

The first ansatz (we consider it the more classical one) is of the form
\begin{equation}\label{eq.ansatz1}
  \sum_{(\alpha,\beta)\in S} c_{\alpha,\beta}M^{\alpha}L^{\beta}
\end{equation}
where the unknown coefficients~$c_{\alpha,\beta}$ have to be determined in~$\BQ(q)$.
We will refer to the finite set~$S\subseteq\BN^2$ as the \emph{structure set} of the ansatz.
It is easy to see that with ansatz~\eqref{eq.ansatz1} at least the first
$(o+1)(d+2)-1$ (resp. $(o+2)(d+2)-2$) values of a sequence are needed in 
order to guess a homogeneous (resp. inhomogeneous) recursion of 
order~$o$ and coefficient degree~$d$, i.e., when $0\leq\alpha\leq d$
and $0\leq\beta\leq o$. Table~\ref{tab.Jdeg} illustrates how fast the
entries of the colored Jones polynomials grow, and
it becomes obvious that we cannot go very far with this ansatz.

The second ansatz is of the form
\begin{equation}\label{eq.ansatz2}
  \sum_{(\alpha,\beta,\gamma)\in S} c_{\alpha,\beta,\gamma}q^{\gamma}M^{\alpha}L^{\beta}
\end{equation}
where the unknowns $c_{\alpha,\beta,\gamma}$ are elements of~$\BQ$
and $S\subseteq\BN^3$ is again a finite structure set.
This alternative is particularly promising when the $q$-degrees
of the sequence grow very fast as it is the case with the
colored Jones polynomials. On the one hand, we have many more
unknowns than in~\eqref{eq.ansatz1}, but already the first
$30$ values of~$J_{p,n}(q)$ suffice to generate thousands of equations.
However, this method has its limits, too, as it can be seen from 
Table~\ref{tab.results}.
Guessing the last entry ($p=5$) would require to solve a linear
system over~$\BQ$ with $17\cdot 289\cdot 2175$ (about 10 million) unknowns.

The small instances in our family of problems can be done with either technique
and without further previous knowledge. But for finding the larger recurrences
presented in this paper, some optimizations are necessary.
What helped considerably in reducing the size of the computations,
is the fact that certain properties of the structure sets in~\eqref{eq.ansatz1}
or~\eqref{eq.ansatz2} can be deduced a priori. In particular, the 
AJ Conjecture of \cite{Ga1}, and the Newton polygon of the $A$-polynomial of 
the pretezel knots $K_p$ from \cite{Ma} allow us to make guess of the 
exponents $(\a,\b)$ that appear in Equation \eqref{eq.ansatz2}, up to
an overall translation in the $M$-direction.
The only missing link to obtain the structure set is how far we have
to translate the Newton polygon in $M$-direction (this has to be found
out by trial and error), see Figure~\ref{fig.newton}. Obviously this
translated Newton polygon has much fewer lattice points than the
rectangular box, and this helps a lot in the computations.

\begin{figure}[htpb]
$$
\psdraw{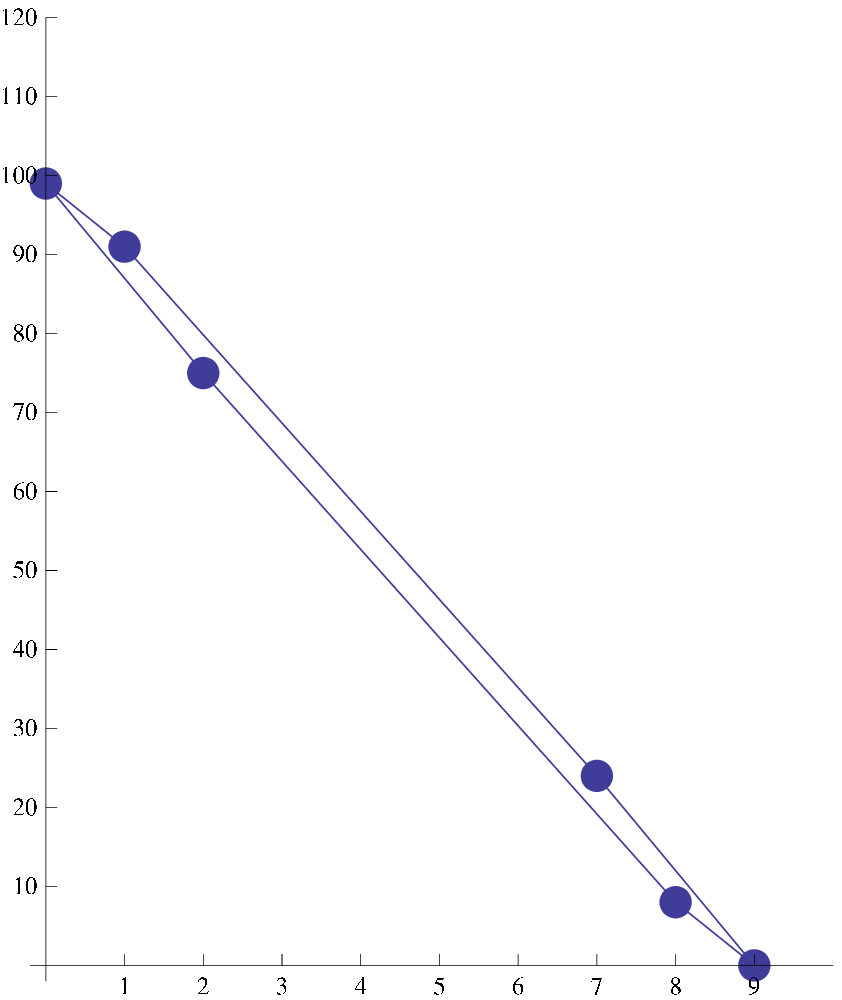}{2.5in}
\psdraw{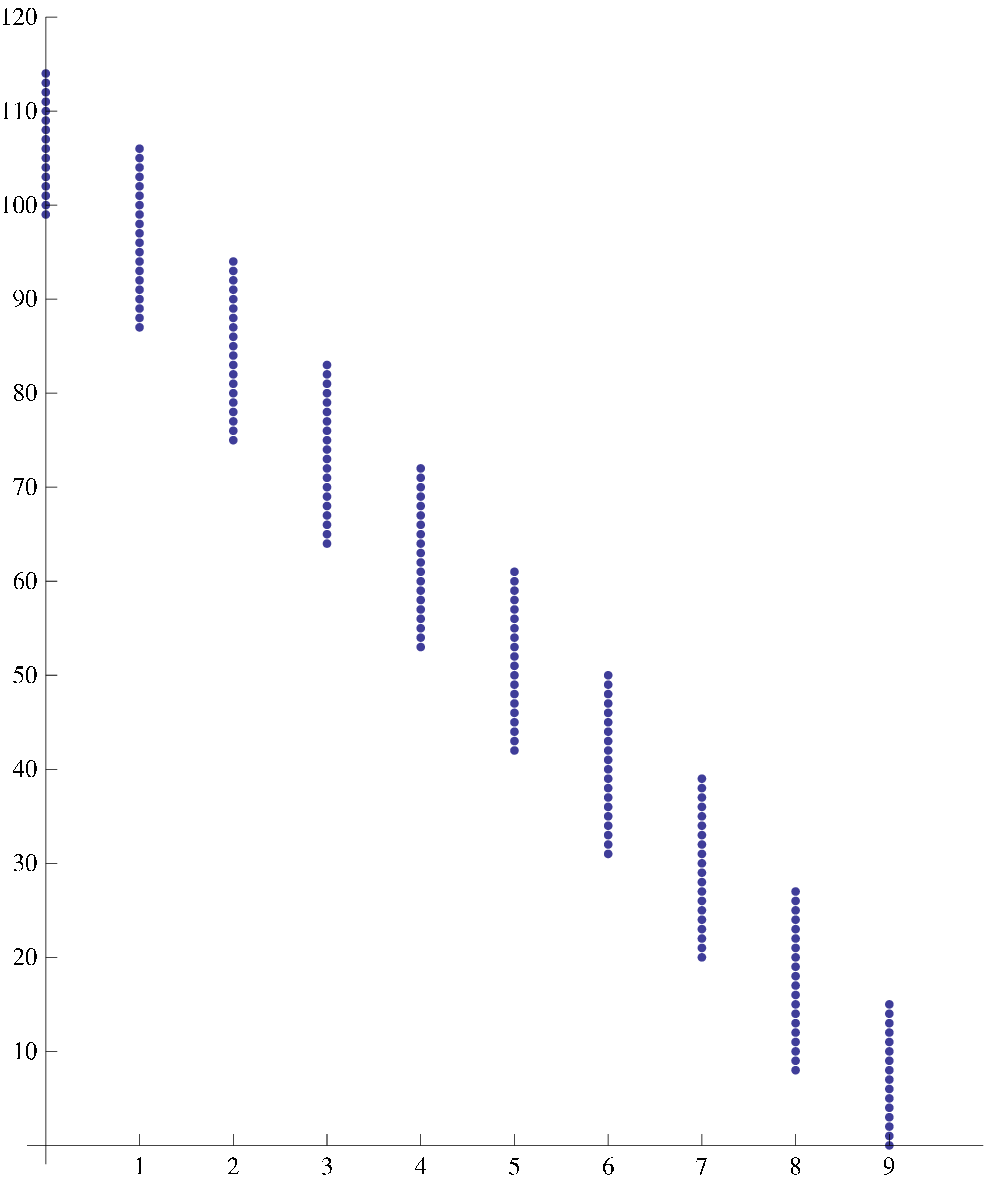}{2.5in}
$$
\caption{The Newton polygon for the (-3,2,9) pretzel knot~$K_3$ (left) 
and the structure set of the recurrence for its colored Jones 
polynomial (right)}
\lbl{fig.newton}
\end{figure}

The second trick that allows us to compute the large recurrences
listed in Table~\ref{tab.results} is to use ansatz~\eqref{eq.ansatz1}
in connection with modular computations. This means that we compute
the sequence~$J_{p,n}(q)$ only for specific integers~$q$ and modulo
some prime number~$m$. Then the monstrous Laurent polynomials (see
Tables~\ref{tab.cjk0} and~\ref{tab.Jdeg}) shrink to a natural number
between~$0$ and~$m-1$. Thus we can compute the first few hundred
values of the colored Jones sequence ``easily'' (i.e., in a few
hours). The guessing then is done efficiently, since only a single
nullspace computation modulo~$m$ of a moderately-sized matrix (less
than 1000 columns) is required.  This procedure has to be performed
for many specific values of~$q$ to be able to interpolate and
reconstruct the polynomial coefficients~$c_{\alpha,\beta}$.  In
general, we also would have to use several prime numbers in order to
recover the integer coefficients via chinese remaindering and rational
reconstruction. However, the integer coefficients in the present
problems are so small (see Table~\ref{tab.results}) that a single
prime suffices to recover the coefficients in~$\BZ$. Note that this
strategy is perfectly suited for parallel computations.

We have seen that the number of sequence entries that we need to
compute is determined by the number of unkowns in the ansatz, and the
number of interpolation points (specific values for~$q$) by the
$q$-degree of the coefficients in the recurrence. We mention some
further tricks for reducing these two parameters. The $q$-degree can
be decreased by considering a slight variation of the original
sequence~$f(n)$, namely $g(n):=f(n+s)$ for some $s\in\BZ$.  Even more,
a similar trick can be used to halve the number of unknowns, by
exploiting the fact the there is an $s=t/2$ for $t\in\BZ$ such that
the substitution $n\to n+s$ in the coefficients of the recurrence
reveals the following symmetry:
\[
  c_{\alpha,\beta}=c_{m-\alpha,l-\beta}\text{ for }t\text{ even,}\qquad\text{and }
  c_{\alpha,\beta}=-c_{m-\alpha,l-\beta}\text{ for }t\text{ odd},
\]
where $m$ is the $M$-degree and $l$ is the $L$-degree of the recurrence.
The above symmetry is equivalent to the palindromic property of
the sought recursion.

We want to remark that most of the computation time was used to compute
the data which then later served as input for the guessing. In the most
difficult example presented here, the recurrence for $K_5$, it took
1607 CPU days to produce the first $744$ entries of $J_{5,n}(q)$
for $700$ different values of~$q$. Since we ran this computation on a
cluster with several hundred processors, it finished within a few days.

\begin{table}[htpb]
{\tiny
\begin{tabular}{|l|c|c|c|c|c|c|c|c|c|c|c|}  \hline
$p$\rule{0em}{1em} & $-5$ & $-4$ & $-3$ & $-2$ & $-1$ & $0$ & $1$ & $2$ & $3$ & $4$ & $5$ \\ \hline
\texttt{ByteCount}\rule{0em}{1em} & $5.7\times 10^7$ & $1.1\times 10^7$ & $1.1\times 10^6$ & $32032$ & 
  $1192$ & $1616$ & $1616$ & $47016$ & $2.3\times 10^6$ & $1.9\times 10^7$ & $8.6\times 10^7$ \\ \hline
$L$-degree\rule{0em}{1em} & $12$ & $9$ & $6$ & $3$ & $1$ & $2$ & $2$ & $6$ & $9$ & $12$ & $15$ \\ \hline
$M$-degree\rule{0em}{1em} & $125$ & $66$ & $27$ & $12$ & $6$ & $13$ & $16$ & $58$ & $114$ & $191$ & $288$ \\ \hline
$q$-degree\rule{0em}{1em} & $946$ & $392$ & $85$ & $19$ & $3$ & $13$ & $16$ & $233$ & $514$ & $1151$ & $2174$ \\ \hline
largest cf.\rule{0em}{1em} & $3.0\times 10^8$ & $12345$ & $33$ & $4$ & $1$ & $2$ & $2$ & $6$ & $118$ & 
  $386444$ & $2.2\times 10^{11}$ \\ \hline
translation\rule{0em}{1em} & $68$ & $39$ & $18$ & $5$ & $1$ & $1$ & $1$ & $3$ & $15$ & $36$ & $65$ \\ \hline
\end{tabular}
}
\caption{Some data concerning the recursion relations for the colored Jones polynomial of the pretzel knots~$K_p$, $p=-5,\dots,5$:
the size of the recurrence (in bytes, using the Mathematica command \texttt{ByteCount}), the order (or $L$-degree), 
the coefficient degree (or $M$-degree), the degree of~$q$ in the coefficients, the largest integer coefficient,
and how far the Newton polygon had to be translated in order to find this recurrence.}
\lbl{tab.results}
\end{table}

\appendix

\section{The recursion for $K_{-2}$ and $K_2$}
\lbl{sec.K2}

The operators $A_{\pm 2}(M,L,q)$ and $b_{\pm 2}(M,q)$ are given as follows.

{\small

\vspace{0.1in}
\noindent \begin{math}
A_{-2}(M,L,q)= 
-(q^2 M-1) (q^2 M+1) (q^4 M-1) (q^3 M^2-1)L^3 
-q (q^3 M-1)^2 (q^3 M+1) (q^3 M^2-1) (q^{14} M^5-q^{11} M^4-(q^{10}-q^9-q^8+q^7)M^3 + (q^7+q^4)M^2+2 q^3 M-1)L^2 
+q^7 M^2 (q^2 M-1)^2 (q^2 M+1) (q^7 M^2-1) (q^{11} M^5-2 q^9 M^4 -(q^8+q^5)M^3 + (q^6-q^5-q^4+q^3)M^2+q^2 M-1)L 
-q^{16} M^7 (q M-1) (q^3 M-1) (q^3 M+1) (q^7 M^2-1) 
\end{math}

\vspace{0.1in}
\noindent \begin{math}
b_{-2}(M,q) = q^6 M^2 (q^2 M+1) (q^3 M+1) (q^3 M^2-1) (q^5 M^2-1) (q^7 M^2-1)
\end{math}

\vspace{0.1in}
\noindent \begin{math}
A_2(M,L,q)= q^{59} (q^2 M-1) (q^3 M-1) (q^7 M-1) L^6 -q^{112} M^8 (q^2 M-1) (q^3 M-1) (q^6 M-1)^3 L^5
-q^{167} M^{18} (q^2 M-1) (q^5 M-1)^2 (q^5 M+q+1) L^4+(q-1) (q+1) q^{207} M^{27} (q^2 M-1) (q^4 M-1)^2 (q^6 M-1) L^3
+q^{240} M^{36} (q^3 M-1)^2 (q^4 M+q^3 M+1) (q^6 M-1) L^2+q^{263} M^{45} (q^2 M-1)^3 (q^5 M-1) (q^6 M-1) L
-q^{279} M^{55} (q M-1) (q^5 M-1) (q^6 M-1)
\end{math}

\vspace{0.1in}
\noindent \begin{math}
b_2(M,q)= q^{89} M^5 (q^{192} M^{48}-q^{186} (q+1) M^{47}+q^{181} M^{46}-q^{187}M^{45}+q^{181} (q+1) M^{44}-q^{176} (q^7+1) M^{43}
 +q^{177} (q^4+q+1) M^{42}-q^{172} (q^7+q^4+q^3+1) M^{41}+q^{170} (q^4+q^3+1)M^{40}+q^{168} (q^6-1) M^{39}
 -q^{168} (-q^3+q+1) M^{38}-q^{163}(q^6+q^3+q^2-1) M^{37}+q^{160} (q^4+q^3+1) M^{36}-q^{158} (q^5+1) M^{35}
 +q^{157} (q^4-q^3+q+1) M^{34}+q^{152} (q^6-q^4+q^2-1)M^{33}-q^{149} (q^3-q+1) M^{32}-q^{148} (q^3+1) M^{31}
 +q^{142} (q^3+1) M^{30}+q^{142} (q^2-1) M^{29}+q^{136} (q^4-q^2+1) M^{28}+q^{134} (q^2-1) M^{27}-q^{134} M^{26}-q^{124} M^{25}
 +2 q^{122} M^{24}-q^{116} M^{23}-q^{118} M^{22}+q^{110} (q^2-1) M^{21}+q^{104}(q^4-q^2+1) M^{20}+q^{102} (q^2-1) M^{19}+
 q^{94} (q^3+1) M^{18}-q^{92} (q^3+1) M^{17}-q^{85} (q^3-q+1) M^{16}+q^{80}(q^6-q^4+q^2-1) M^{15}+q^{77} (q^4-q^3+q+1) M^{14}
 -q^{70} (q^5+1)M^{13}+q^{64} (q^4+q^3+1) M^{12}-q^{59} (q^6+q^3+q^2-1) M^{11}-q^{56} (-q^3+q+1) M^{10}+q^{48} (q^6-1) M^9
 +q^{42} (q^4+q^3+1) M^8-q^{36} (q^7+q^4+q^3+1) M^7+q^{33} (q^4+q+1) M^6-q^{24} (q^7+1) M^5+q^{21} (q+1) M^4
 -q^{19} M^3+q^5 M^2-q^2 (q+1) M+1)
\end{math}
}

\vspace{0.1in}
\noindent The values of $(A_p(M,L,q),b_p(M,q))$ for $p=-5,\dots,5$ are 
available from
\begin{center}
{\tt \url{http://www.risc.jku.at/people/ckoutsch/pretzel/}}
\end{center}
or
\begin{center}
{\tt \url{http://www.math.gatech.edu/~stavros/publications/pretzel.data/}}
\end{center}

\section{The $M$-factors}
\lbl{sec.Mfactors}

The $M$-factors $\e_p(M) \in \BQ[M]$ given as follows for $p=-5,\dots,5$.

{\small
\vspace{0.1in}
\noindent \begin{math}
\e_{-5}(M)=
-(-1+M)^9 (1+M)^5 (1+M+M^2) (1-2 M+5 M^2+6 M^3-14 M^4+20 M^5+17 M^6-48 M^7+43 M^8
+40 M^9-67 M^{10}+40 M^{11}+43 M^{12}-48 M^{13}+17 M^{14}+20 M^{15}-14 M^{16}+6 M^{17}+5 M^{18}-2 M^{19}
+M^{20}) (12-32 M+34 M^2+18 M^3-171 M^4+462 M^5-680 M^6+240 M^7+1054 M^8-2126 M^9
+1332 M^{10}+1080 M^{11}-3016 M^{12}+2558 M^{13}-227 M^{14}-2256 M^{15}+3187 M^{16}-2256 M^{17}-227 M^{18}
+2558 M^{19}-3016 M^{20}+1080 M^{21}+1332 M^{22}-2126 M^{23}+1054 M^{24}+240 M^{25}-680 M^{26}+462 M^{27}
-171 M^{28}+18 M^{29}+34 M^{30}-32 M^{31}+12 M^{32})
\end{math}

\vspace{0.1in}
\noindent \begin{math}
\e_{-4}(M)=
(-1+M)^9 (1+M)^6 (1+M^4)^2 (4-2 M+7 M^2+10 M^3-14 M^4+34 M^5-7 M^6+6 M^7+24 M^8+6 M^9-7 M^{10}+34 M^{11}-14 M^{12}+10 M^{13}+7 M^{14}-2 M^{15}+4 M^{16})
\end{math}

\vspace{0.1in}
\noindent \begin{math}
\e_{-3}(M)=
2 (-1+M)^5 (1+M)^5 (1-M+M^2)^3 (1+M+M^2)
\end{math}

\vspace{0.1in}
\noindent \begin{math}
\e_{-2}(M)=
-(-1+M)^3 (1+M)^2
\end{math}

\vspace{0.1in}
\noindent \begin{math}
\e_{-1}(M)=
-1+M
\end{math}

\vspace{0.1in}
\noindent \begin{math}
\e_{0}(M)=
-1+M
\end{math}

\vspace{0.1in}
\noindent \begin{math}
\e_{1}(M)=
1-M
\end{math}

\vspace{0.1in}
\noindent \begin{math}
\e_{2}(M)=
(-1+M)^3
\end{math}

\vspace{0.1in}
\noindent \begin{math}
\e_{3}(M)=
2 (-1+M)^5 (1+M)^2 (1+M+M^2)^2 (1-6 M+13 M^2-6 M^3+M^4)
\end{math}

\vspace{0.1in}
\noindent \begin{math}
\e_{4}(M)=
-(-1+M)^7 (1+M)^3 (1+M^2) (1-4 M+4 M^2+4 M^3-8 M^4+4 M^5+4 M^6-4 M^7+M^8) (4-14 M
+33 M^2+12 M^3-330 M^4+328 M^5-9 M^6+226 M^7+650 M^8+226 M^9-9 M^{10}+328 M^{11}-330 M^{12}
+12 M^{13}+33 M^{14}-14 M^{15}+4 M^{16})
\end{math}

\vspace{0.1in}
\noindent \begin{math}
\e_{5}(M)=
(-1+M)^9 (1+M)^2 (1+M+M^2) (1-2 M+4 M^2-2 M^3-M^4-109 M^6-78 M^7+406 M^8+162 M^9
-417 M^{10}+162 M^{11}+406 M^{12}-78 M^{13}-109 M^{14}-M^{16}-2 M^{17}+4 M^{18}-2 M^{19}+M^{20}) (12-88 M
+318 M^2-698 M^3+381 M^4+4924 M^5-20623 M^6+30440 M^7+4694 M^8-61074 M^9+47268 M^{10}
+15096 M^{11}-11350 M^{12}-42702 M^{13}+37078 M^{14}+55502 M^{15}-112131 M^{16}+55502 M^{17}+37078 M^{18}
-42702 M^{19}-11350 M^{20}+15096 M^{21}+47268 M^{22}-61074 M^{23}+4694 M^{24}+30440 M^{25}-20623 M^{26}
+4924 M^{27}+381 M^{28}-698 M^{29}+318 M^{30}-88 M^{31}+12 M^{32})
\end{math}
}

\vspace{0.1in}
An explicit calculation shows that $\e_p(M)$ are palindromic polynomials,
i.e., $\e_p(M)/e_p(1/M)=-M^{\d_p}$ where $\d_p$ is given by
$$
\{68, 39, 18, 5, 1, 1, 1, 3, 15, 36, 65\}
$$
for $p=-5,\dots,5$.
One factor of $\e_p(M)$ is easy to spot, namely it is the Alexander polynomial
$\D_p(M)$, in accordance with the loop expansion; see Equation \eqref{eq.loop}.
The other factors of $\e_p(M)$ are more mysterious with no clear geometric
definition. Note finally that $\d_p$ given above agree with the translation
factor of the Newton polygon of $K_p$ given in the last row of Table 
\ref{tab.results}. This agreement is consistent with the fact that our 
computed recursion relations $(A_p(M,L,q),b_p(M,q))$ are palindromic
for $p= -5, \dots, 5$.

\subsection*{Acknowledgment}
The paper came into maturity following a request by D. Zagier for an
explicit formula for the recursion of the colored Jones polynomial of
a knot, during visits of the first author in the Max-Planck-Institut
f\"ur Mathematik in 2008-2010. The authors met during the conference
in honor of Doron Zeilberger's 60th birthday at Rutgers University in
2010, and wish to thank the organizers of the Zeilberger Fest, and
especially D. Zagier and D. Zeilberger for their interest and for many
stimulating conversations.

\tableofcontents

\bibliographystyle{hamsalpha}\bibliography{biblio}
\end{document}